\documentclass[12pt,a4paper]{amsart}

\usepackage{amsfonts, amsmath, amssymb, amsthm, amscd, hyperref}

\usepackage{a4wide}

\newtheorem{thm}{Theorem}
\newtheorem{lem}{Lemma}
\newtheorem{con}{Conjecture}

\theoremstyle{definition}
\newtheorem{defn}{Definition}

\theoremstyle{remark}

\newcommand{\pt}{\mathrm{pt}}

\renewcommand{\int}{\mathop{\rm int}}

\renewcommand{\epsilon}{\varepsilon}

\begin{document}

\title[Knaster's problem for $(Z_2)^k$-symmetric subsets\dots]{Knaster's problem for $(Z_2)^k$-symmetric\\ subsets of the sphere $S^{2^k-1}$}

\author{R.N.~Karasev}
\thanks{This research was partially supported by the Dynasty Foundation.}

\email{r\_n\_karasev@mail.ru}
\address{
Roman Karasev, Dept. of Mathematics, Moscow Institute of Physics
and Technology, Institutskiy per. 9, Dolgoprudny, Russia 141700}


\begin{abstract}
We prove a Knaster-type result for orbits of the group $(Z_2)^k$ in $S^{2^k-1}$, calculating the Euler class obstruction. Among the consequences are: a result about inscribing skew crosspolytopes in hypersurfaces in $\mathbb R^{2^k}$, and a result about equipartition of a measures in $\mathbb R^{2^k}$ by $(Z_2)^{k+1}$-symmetric convex fans.
\end{abstract}

\maketitle

\section{Introduction}

Recall the well-known Knaster problem~\cite{kna1947}.

\begin{con}[Knaster's problem for functions]
\label{knaster}
Let $S^{d-1}$ be a unit sphere in $\mathbb R^d$. Suppose we are given $d$ points $x_1, \ldots, x_d\in S^{d-1}$ and a continuous function $f: S^{d-1}\to \mathbb R$. Then there exists a rotation $\rho\in SO(d)$ such that
$$
f(\rho(x_1)) = f(\rho(x_2)) = \dots = f(\rho(x_d)).
$$
\end{con}

The Knaster conjecture turned out to be false in general. Surprisingly, it took quite a long time to find a counterexample. First counterexamples were found in~\cite{mak1986,babo1990,chen1990} to the generalized Knaster conjecture for maps to $\mathbb R^m$ for $m>1$. The counterexamples to Conjecture~\ref{knaster} were found in~\cite{kasha2003,hiri2005} for large enough $d$, some particular sets $\{x_1,\ldots, x_d\}$, and particular functions $f$.

Still, some partial solutions for this problem were given by several authors. For example, in~\cite{dys1951,floyd1955,gri1991,hms2000} the case of $\mathbb R^3$ and some its strengthened variants were solved. In~\cite{yayu1950} the case when $\{x_1,\ldots, x_d\}$ is an orthonormal base was solved positively, in~\cite{mak1990} the case of $Z_p$-symmetric set $\{x_1, \ldots, x_p\}$ on a circle was considered, and in~\cite{vol1992} the Knaster problem was solved for an orbit of a certain action of $(Z_p)^k$ on $\mathbb R^{p^k}$ for odd $p$.

Using the result of~\cite{vol1992} on the Knaster problem, some corollaries on inscribing $(Z_p)^k$-symmetric crosspolytopes and partitioning a measure by $(Z_p)^k$-symmetric ``fans'' of cones were obtained by the author in~\cite{kar2009ins,kar2009mp} for odd primes $p$.

In the present paper the case of $(Z_2)^k$-symmetry is considered. In fact, instead of the group of symmetry $G=(Z_2)^k$, we have to consider some larger group $G'$, such that $G$ is its subgroup and a quotient group (see the proof below). The function $f:S^{d-1}\to\mathbb R$ is assumed to be even, permitting the action of the larger group $G'$ on the corresponding representation.

\begin{thm}
\label{knaster-pow2}
Denote $G=(Z_2)^k$, $k\ge 2$, $q=2^k$, and consider some subgroup $H\subset G$, which is isomorphic to $(Z_2)^{k-1}$. Denote some element $\sigma\in G\setminus H$. 

Suppose that $\mathbb R[G]$ has some $G$-invariant scalar product $(\cdot, \cdot)$ and let $S\subset\mathbb R[G]$ be the unit sphere of this scalar product. Consider a point $x\in S$ such that the set $Gx$ is a base of $\mathbb R[G]$ and $Hx\perp H\sigma x$.

Then for any even continuous function $f : S\to \mathbb R$ there exists a rotation (w.r.t. the product $(\cdot,\cdot)$) $\rho\in SO(q)$ such that the function $f$ is constant on $\rho(Gx)$.
\end{thm}

This result is different from the original formulation of the Knaster's problem, because the function $f$ is considered even. In the case when the vectors $Gx$ form an orthogonal base this restriction can be avoided, and we obtain another proof of Knaster's conjecture for the case of orthogonal bases (see~\cite{yayu1950}). Of course, in~\cite{yayu1950} the dimension was arbitrary, not necessarily a power of two, but the proof of Theorem~\ref{knaster-pow2} is still of value because it leads to some corollaries, presented in the next section. 

It is also possible to prove another version of this theorem.

\begin{thm}
\label{knaster-pow2-odd}
Theorem~\ref{knaster-pow2} is true for an odd function $f$ instead of even.
\end{thm}

The proofs of Theorems~\ref{knaster-pow2} and \ref{knaster-pow2-odd} do not extend to the case when $f$ is neither even nor odd. However it is possible that more detailed computations in cohomology may prove the general case.

\section{Related results}

Knaster's conjecture is a close relative of many other combinatorial geometric problems, including the problem of inscribing polytopes into hypersurfaces and the problem of partitioning mass distributions
by convex fans.

In~\cite{shn1944} a theorem on inscribing a square into a closed smooth curve in $\mathbb R^2$ without self-intersections (the Schnirelmann theorem) was proved. In~\cite{gug1965} the possible generalization of the Schnirelmann theorem was considered: for any image $M$ of a smooth embedding $S^{d-1}\to\mathbb R^d$ there is a regular crosspolytope $C$ (a set, similar to the convex hull of vectors $\pm e_1, \ldots, \pm e_d$, where $(e_1,\ldots, e_d)$ is a orthonormal base) such that, all its vertices lie on $M$. Still, the sketch of the proof in~\cite{gug1965} was incorrect, the main idea was that the inscribed crosspolytope depends continuously on deformations of $M$. But it is not true already in the case of the plane, where the inscribed squares may disappear under the deformations of the curve. 

The difference between the case of the plane and $\mathbb R^d$ for $d>2$ is as follows. In the plane only even number of squares may disappear at a time, while deforming the curve, thus keeping the number of inscribed squares odd. In the space there is a variety of inscribed crosspolytopes, and the counting and parity argument does not work.

In the book~\cite{kleewa1996} the problem on inscribing a regular octahedron ($3$-dimensional crosspolytope) was stated as unsolved. The $3$-dimensional case was solved in~\cite{mak2003}, and the case of dimension $p^k$, where $p$ is an odd prime, was solved in~\cite{kar2009ins}.  Here we consider the case of dimension $q=2^k$. 

Note that the problem of inscribing a crosspolytope to the boundary of a non-smooth convex body does not follow directly from the smooth case. In~\cite{kar2009ins} some restrictions on a non-smooth convex body (``non-acuteness'') were considered in the case of the regular crosspolytope. Those conditions guaranteed the existence of an inscribed regular crosspolytope by approximating the non-smooth body by smooth bodies and going to the limit. We formulate the smooth case here, noting that for the regular crosspolytope and convex bodies the smoothness may be replaced with the same ``non-acuteness'' condition as in~\cite{kar2009ins}.

\begin{thm}
\label{insrccross-pow2}
Let $G$, $H$, $q$, $x$ be the same as in Theorem~\ref{knaster-pow2}. Suppose the hypersurface $M\subset \mathbb R[G]$ is an image of a smooth embedding $S^{q-1}\to\mathbb R[G]$. Then there exists a similarity (a composition of a rigid motion and a homothety) $\rho$ that preserves the orientation, such that $2q$ points $\rho(\pm Gx)$ lie on $M$.
\end{thm}

In particular, this theorem gives a theorem about inscribing regular crosspolytopes in $\mathbb R^q$ for $q=2^k$, because in this way the vectors $(e_1, \ldots, e_d)$ are permuted by some $G=(Z_2)^k$, and the orthogonality condition $Hx\perp H\sigma x$ obviously holds.

Theorems on partitioning a measure or a set of measures in $\mathbb R^d$ into equal parts have quite a long history, see, for example, \cite{st1942, ste1945, grun1960}, in particular, some topological methods were used to give partitions by ``fans'', that is a set of hyperplanes or half-hyperplanes, having a common plane of codimension $2$, see~\cite{mak1994,barmat2001,zivvre2001,ziv2004}. 

Here we consider a partition of $\mathbb R^d$ by a system of cones with apex at the origin, let us call such a system ``fan'' again. Then this system may be allowed to be moved by rigid motions, so as to give a needed partition of a measure. One case of such results was proved for $\mathbb R^3$ in~\cite{mak1988}, a generalization for dimensions $p^k$ for odd prime $p$ was made in~\cite{kar2009mp}, and here we consider the case of dimension $q=2^k$.

\begin{thm}
\label{mespartcones-pow2}
Let $G$, $H$, $q$ be the same as in Theorem~\ref{knaster-pow2}. Suppose that $C$ is a closed cone with apex at the origin, the family of cones $\{\pm g(C)\}_{g\in G}$ gives a partition of $\mathbb R[G]$, the subfamily $\{g(C)\}_{g\in G}$ have a unique common ray. 

Consider the involution $\tau: \mathbb R[G]\to \mathbb R[G]$, defined on the base $Gx$ by 
$$
\forall g\in H,\ \tau(gx) = gx,\quad,\forall g\in G\setminus H,\ \tau(gx) = -gx.
$$
Suppose that $\tau(g(C)) = g(C)$ for $g\in H$ and $\tau(g(C)) = -g(C)$ for $g\in G\setminus H$.

Then for any absolutely continuous probabilistic measure $\mu$ on $\mathbb R[G]$ there is a rigid motion $\rho$ that preserves the orientation, such that for any $g\in G$
$$
\mu(\rho(g(C))) = \mu(\rho(-g(C)) = \frac{1}{2d}.
$$
\end{thm}

For example, in this theorem we can take the Voronoi cells of the system $\{\pm gx\}_{g\in G}$ from Theorem~\ref{knaster-pow2} as the cones.

The case $k = 1$ ($q=2$, $\mathbb R[G]$ is the plane) is not considered here, because in this case Theorem~\ref{insrccross-pow2} is the Schnirelmann theorem~\cite{shn1944}. The plane case of Theorem~\ref{mespartcones-pow2} is also well-known, e.g. see~\cite{grun1960}, Remark~4.v. 

Nevertheless, the proof of Theorem~\ref{knaster-pow2} can be modified for the case $q=2$. In this case the group $G$ and the involution $\tau$ alter the orientation of the plane and we have to consider $O(2)$ instead of $SO(2)$ in all the reasonings.  

\section{Equivariant cohomology of $G$-spaces}
\label{eq-cohomology}

We consider topological spaces with continuous (left) action of a finite group $G$ and continuous maps between such spaces that commute with the action of $G$. We call them $G$-spaces and $G$-maps.

For basic facts about (equivariant) topology and vector bundles the reader is referred to the books~\cite{milsta1974,hsiang1975,mishch1998}. The cohomology is taken with coefficients in $Z_2$, in notations we omit the coefficients. Let us start from some standard definitions.

\begin{defn}
Denote $EG$ the classifying $G$-space, which can be thought of as an infinite join $EG=G*\dots *G*\dots$ with diagonal left $G$-action. Denote $BG=EG/G$. For any $G$-space $X$ denote $X_G=(X\times EG)/G$, and put (\emph{equivariant cohomology in the sense of Borel}) $H_G^*(X) = H^*(X_G)$. It is clear that for a free $G$-space $X$ the space $X_G$ is homotopy equivalent to $X/G$. 
\end{defn}

Consider the algebra of $G$-equivariant cohomology of the point $A_G = H_G^*(\pt) = H^*(BG)$. For any $G$-space $X$ the natural projection $\pi_X : X\to\pt$ induces the natural map of cohomology $\pi_X^* : A_G\to H_G^*(X)$.

For a group $G=(Z_2)^k$ the algebra $A_G$ (see~\cite{hsiang1975}) is the algebra of polynomials of $k$ one-dimensional generators $v_i$.

We are going to find the equivariant cohomology of a $G$-space $X$ using the following spectral sequence (see~\cite{hsiang1975,mcc2001}).

\begin{thm}
\label{specseqeq}
The natural fiber bundle $\pi_{X_G} : X_G\to BG$ with fiber $X$ gives the spectral sequence with $E_2$-term
$$
E_2^{x, y} = H^x(BG, \mathcal H^y(X)),
$$
that converges to the graded module, associated with the filtration of $H_G^*(X)$.

The system of coefficients $\mathcal H^y(X)$ is obtained from the cohomology $H^y(X)$ by the action of $G = \pi_1(BG)$. The differentials of this spectral sequence are homomorphisms of $H^*(BG)$-modules.
\end{thm}

We need the following lemma, stated in~\cite{vol1992} based on results from~\cite{dic1911,mui1975}.

\begin{lem}
\label{euler-nz}
Let $G=(Z_2)^k$, $q=2^k$, and let $I[G]$ be the subspace of $\mathbb R[G]$, consisting of elements
$$
\sum_{g\in G} a_g g,\quad \sum_{g\in G} a_g = 0.
$$
Then the only nonzero Stiefel-Whitney classes are $w_{q-2^l}(I[G])\in A_G$ $(l=0,\ldots, k)$, and $w_{q-1}(I[G])$ is not contained in the ideal of $A_G$, generated by $w_k(I[G])$ with $k<q-1$. 
\end{lem}

\section{Proof of Theorem~\ref{knaster-pow2}}

Consider the involution $\tau$, as it is defined in the statement of Theorem~\ref{mespartcones-pow2}. Since $Hx\perp H\sigma x$, $\tau$ acts orthogonally on $\mathbb R[G]$. The group $G$ and the involution $\tau$ generate a (non-commutative) group $G'$ of order $2^{k+2}$. 

Denote the natural inclusion $\iota : G\to G'$. It is easy to see that there is a homomorphism $\pi : G'\to G$ with kernel $\{1,\tau,-\tau,-1\}$, left inverse to $\iota$. This homomorphism is given by considering the permutation by $G'$ of $1$-dimensional subspaces, corresponding to vectors $gx$ ($g\in G$). 

As it was done in the approach to the Knaster conjecture in~\cite{vol1992}, we consider the map $\phi : SO(q)\to \mathbb R[G]$ defined as follows:
$$
\phi(\rho) = \sum_{g\in G} f(\rho(gx)) g. 
$$
If $G$ acts on $SO(q)$ by right multiplications by $g^{-1}$ (this is a left action) and acts by left multiplications on $\mathbb R[G]$, then $\phi$ is a $G$-map. We need to prove that $\phi$ maps some $\rho$ to an element $\sum_{g\in G} cg$ ($c$ is a constant).
Taking the orthogonal projection, we obtain a map $\psi : SO(q)\to I[G]$, and for this map we only need to prove that it maps some $\rho$ to zero.

Similar to the considerations for the case of odd $p$ and the group $(Z_p)^k$, it would be natural to prove that the Stiefel-Whitney class $w_{q-1}(I[G])\not=0\in H_G^{q-1}(SO(q))$. In fact, this is not true, the class is zero, and actually we have to consider the group $G'$ and prove that $w_{q-1}(I[G])\not=0\in H_{G'}^{q-1}(SO(q))$. Here we consider $I[G]$ as a representation of $G'$ given by the map $\pi : G'\to G$, so $\phi$ and $\psi$ become $G'$-maps. This is possible since the function $f$ is even, $\pm 1$ and $\pm \tau$ either stabilize $gx$ (for $g\in G$) or map it to its opposite $-gx$, hence
$$
\phi(\rho(\pm\tau)) = \sum_{g\in G} f(\rho(\pm\tau gx)) g = \sum_{g\in G} f(\rho(\pm gx)) g = \sum_{g\in G} f(\rho(gx)). 
$$

We also consider a subgroup $G''\subset G'$, generated by $\tau$ and $H\subset G$. This group is isomorphic to $(Z_2)^k$. In order to distinguish the action of all these groups on the original space $\mathbb R[G]$ and the other actions we fix some notation.

\begin{defn}
\label{actions}
Denote the original space $\mathbb R[G]$ with $G'$ action by $U$. Denote $I[G]$ with the natural $G$-action by $V$. Denote the same space with the action of $G'$, induced from the map $\pi : G'\to G$ by $V'$. Denote the space $V'$ with the action of $G''$ induced by the inclusion $\iota'': G''\to G'$ by $V''$.
\end{defn} 

By Lemma~\ref{euler-nz} we have $w_{q-1}(V)\not=0\in A_G$. The action of $G''$ on $V''$ has a nontrivial fixed vector $\sum_{g\in H} g - \sum_{g\in G\setminus H} g$ and therefore $w_{q-1}(V'')=0\in A_{G''}$.

Since there is an inclusion $\iota: G\to G'$, right inverse to the projection $\pi:G'\to G$, then $A_G$ is a summand in $A_{G'}$, so $w_{q-1}(V')\not=0\in A_{G'}$. 

Now we have to prove that $w_{q-1}(V')$ survives in the spectral sequence of Theorem~\ref{specseqeq} for the group $G'$ acting on $SO(q)$. In this case the class $w_{q-1}(\pi_{SO(q)}^*(V'))\not=0\in H_{G'}^{q-1}(SO(q))$ and the map $\psi$ must have a zero.

Note that the groups $G'', G', G$ act trivially on the cohomology $H^*(SO(q))$, because their action can be obviously extended to the action of $SO(q)$ on itself, and $SO(q)$ is connected. Thus the spectral sequence is multiplicative.

Note that the map $\mathfrak I : SO(q)\to SO(q)$ defined by $\mathfrak I(x) = x^{-1}$ turns the considered $G$-action on $SO(q)$ into the standard left multiplication. Denote $SO(q)$ with $G$-action by left multiplication by $SO(q)'$ to distinguish between the two cases. Obviously, the map $\mathfrak I$ gives a natural isomorphism between $G$-spaces $SO(q)$ and $SO(q)'$, and between the spectral sequences of Theorem~\ref{specseqeq}.
 
Consider the natural commutative diagram of maps
$$
\begin{CD}
\label{leseq}
SO(q)'_G @>>> ESO(q)\\
@V{\pi}VV @V{\pi}VV\\
BG @>>> BSO(q),
\end{CD}
$$
where the horizontal arrows are induced by the inclusion $G\to SO(q)$, and the vertical arrows are the fiber bundle projections of Theorem~\ref{specseqeq}. All differentials of the spectral sequence of Theorem~\ref{specseqeq} for the fiber bundle $ESO(q)\to BSO(q)$ are generated by transgressions that kill the Stiefel-Whitney classes $w_m(SO(q))\in H^*(BSO(q))$ ($m=2,\ldots, q$), see Proposition~23.1 in~\cite{bor1953}. Similar to what is done in~\cite{vol1992} for the cohomology with coefficients in $Z_p$ for $p>2$, we conclude that all differentials of the spectral sequence of Theorem~\ref{specseqeq} for $G$-action on $SO(q)'$ (and therefore on $SO(q)$) are generated by transgressions that send the primitive (in terms of~\cite{bor1953}) generators of $H^*(SO(q))$ to the Stiefel-Whitney classes $w_m(SO(q))=w_m(U)\in A_{G}$. The same is true for $G'$ instead of $G$.

The map $\pi : G'\to G$ induces a map between the spectral sequences $\pi^* : E_*^{*,*}(G, SO(q))\to E_*^{*,*}(G', SO(q))$. Obviously, in $E_*^{*,*}(G, SO(q))$ we have the equality $w_{q-1}(SO(q))=w_{q-1}(V)$, and $w_{q-1}(V)$ is ``killed'' by the $q-1$-th transgression. Lemma~\ref{euler-nz} along with the natural projection $\iota^*: A_{G'}\to A_G$ shows that the class $w_{q-1}(V')$ is not contained in the ideal of $A_{G'}$, generated by $w_k(V')$ with $k<q-1$. Thus in the spectral sequence $E_*^{*,*}(G', SO(q))$ the class $w_{q-1}(V')$ survives until the $q-1$-th transgression in $E_*^{*,*}(G', SO(q))$, but still can be ``killed'' by $d_{q-1}$.

Note that there is no more than one primitive generator of $H^*(SO(q))$ in each dimension, hence the differential $d_{q-1}$ in $E_{q-1}^{*,*}(G', SO(q))$ has the image of dimension $1$ in $E_{q-1}^{q-1, 0}(G', SO(q))$. We are going to show that its image does not equal $w_{q-1}(V')$.

Assume the contrary: $d_{q-1} (E_{q-1}^{0, q-2}(G', SO(q))) = \{0, w_{q-1}(V')\}$, since it is a one-dimensional $Z_2$-space. Consider the natural map between the spectral sequences $\iota''^* : E_*^{*,*}(G', SO(q))\to E_*^{*,*}(G'', SO(q))$. It was already mentioned that $\iota''^*(w_{q-1}(V')) = w_{q-1}(V'') = 0$, so from the naturality of $d_{q-1}$ we have
$$
d_{q-1} (E_{q-1}^{0, q-2}(G'', SO(q))) = \{0\}.
$$
But the action of $G''$ on $U$ is isomorphic to the action of $G''$ on $\mathbb R[G'']$. To show this it is sufficient to change the base in $U$ from $\{gx\}_{g\in G}$ to $\{hx + h\sigma x, hx - h\sigma x\}_{h\in H}$ and note that $\tau$ permutes the vectors $hx + h\sigma x$ and $hx - h\sigma x$. So $w_{q-1}(SO(q))=w_{q-1}(\mathbb R[G''])\not=0\in A_{G''}$, and $w_{q-1}(SO(q))\in A_{G''}$ is in the non-trivial image of $d_{q-1}$ of $E_{q-1}^{*,*}(G'', SO(q))$ by Lemma~\ref{euler-nz}. We have come to the contradiction.

\section{Proof of Theorem~\ref{knaster-pow2-odd}}

Consider again the map $\phi : SO(q)\to \mathbb R[G]$ defined as follows:
$$
\phi(\rho) = \sum_{g\in G} f(\rho(gx)) g. 
$$

In the case of odd $f$ this map can be considered as a $G'$-map, the action of $G'$ on the target space $\mathbb R[G]$ becomes isomorphic to its action on $\mathbb R[G]$ from the statement of the theorem, i.e. $\tau g=g$ if $g\in G\setminus H$ and $\tau g = -g$ if $g\in H$. Denote both $G'$-spaces by $U$. 

Consider a subgroup $F=C\times G\subset G'$, where $C=\{1, -1\}$. Its cohomology by the K\"unneth formula is $A_F = A_C\otimes A_G$. We identify $A_C=Z_2[\alpha]$ and consider $A_C=A_C\otimes 1$ and $A_G=1\otimes A_G$ as subalgebras of $A_F$. We are going to find the Stiefel-Whitney classes of $U$ in $A_F$.

\begin{lem}
\label{euler-nz-o}
Denote $q=2^k$. Then the classes $w_{q-2^l}(U)=w_{q-2^l}(\mathbb R[G])\not=0\in A_G\subset A_F$ $(l=0,\ldots, k)$, 
$$
w_q(U) = \sum_{l=0}^k \alpha^{2^l} w_{q-2^l}(U),
$$
and the other classes are zero. The class $w_q(U)$ is (obviously) not contained in the ideal, generated by $w_k(U)$ with $k<q$.
\end{lem}

The proof is based on the fact that $U=L\otimes\mathbb R[G]$, where $L$ has nontrivial $C$-action, and the direct computations of $w(\mathbb R[G])$ in~\cite{dic1911,mui1975,vol1992}, already mentioned in Lemma~\ref{euler-nz}.

Lemma~\ref{euler-nz-o} shows that if we consider the map $\phi$ as a section of a $G'$-bundle, then its Euler class is zero. If $\phi$ still maps some $\rho$ to zero, the proof is complete. 

Assume the contrary: the map $\phi$ does not map any point to zero. Thus it can be projected to the sphere of $U$ to give the $G'$-map $\tilde\phi : SO(q)\to S(U)$. If we consider the map of spectral sequences $\tilde\phi^*: E_*^{*,*} (F, S(U))\to E_*^{*,*} (F, SO(q))$ we note that the Euler class $w_q(U)$ is killed in both spectral sequences by $d_q$. Since in $E_*^{*,*} (F, S(U))$ this is the only nontrivial differential that sends the generator of $H^{q-1}(S(U))$ to $w_q(U)$, we deduce that the map $\tilde\phi^*:H^{q-1}(S(U))\to H^{q-1}(SO(q))$ must be nontrivial.

If the map $\tilde\phi^*:H^{q-1}(S(U))\to H^{q-1}(SO(q))$ is nontrivial, the map $\tilde\phi$ itself must be surjective. It means, in particular, that for some $\rho$ we have
$$
\phi(\rho) = \sum_{g\in G} cg,\quad c\in\mathbb R,
$$
and the proof is complete.

\section{Proofs of the corollaries}

Note that the group $G'$, defined in the previous section, acts naturally on the points $\pm Gx$ and on the cones $\pm gC$ from the respective theorems.

The proofs of these theorems follow the proofs in~\cite{kar2009ins,kar2009mp}. Here we give a sketch of the proof for Theorem~\ref{mespartcones-pow2}, noting the difference from the proof of Theorem~1 in~\cite{kar2009mp}.

A rigid motion $\rho$ of $\mathbb R^d$, that preserves the orientation, has the form
$$
\rho(x) = s(x) + t,
$$
where $s\in SO(d)$ is a rotation and $t\in\mathbb R^d$ is the translation vector. Thus, topologically $E=SO(d)\times \mathbb R^d$. 

Now we define a continuous map $\phi$ from $E$ to $V=\mathbb R[G]\oplus\mathbb R[G]$. For $\rho\in E$ and any $g\in G$ put
$$
a_g = \mu(\rho(g(C)),\quad b_g = \mu(\rho(-g(C)),
$$
$(a_g)$ and $(b_g)$ being the coordinates in $V=\mathbb R[G]\oplus\mathbb R[G]$. The measure $\mu$ is absolutely continuous, so the map $\phi$ is continuous. Note that $G'$ (the same group as in the proof of Theorem~\ref{knaster-pow2}) acts on $E$ by right multiplications by $g^{-1}$. The action of $G'$ on $V$ is the following: $G$ acts by left multiplications and $\tau$ permutes $a_g$ and $b_g$ when $g\not\in H$ and keeps $a_g$ and $b_g$ when $g\in H$. It can be easily seen that the map $g$ is $G'$-equivariant (commutes with the $G'$-action).

Take new coordinates in $V$ as $s_g = a_g + b_g, t_g = a_g - b_g$. Let us enumerate the elements of $G$ some way as $\{g_1, \ldots, g_d\}$. Consider the one-dimensional subspace $L\subset V$, given by equalities
$$
t_{g_1}=\dots=t_{g_d} = 0,\quad s_{g_1}=\dots=s_{g_d}.
$$
In the space $V/L$ denote the $d$-dimensional linear hull of $\{t_{g_1},\ldots, t_{g_d}\}$ by $U$, and the $d-1$-dimensional linear hull of $\{s_{g_1},\ldots, s_{g_d}\}$ by $W$. Let us denote the natural projection $\pi : V\to V/L$ and $f = \pi\circ \phi$. Now all we need is to prove that the map $f: E\to U\oplus W$ maps some motion $\rho$ to zero. 

Let us take a closer look at the action of $G'$ on $U$ and $W$. They are $\tau$-invariant, $\tau$ has some nontrivial action on $U$, while it acts trivially on $W$. Thus the representation $W$ of $G'$ is the same as $V'$ in Definition~\ref{actions}.

The map $f$ can be regarded as a section of the trivial (but not $G'$-trivial) $G'$-bundle ever the $G'$-space. The first obstruction to the existence of a nonzero section is the Euler class of the bundle, in fact it is the highest Stiefel-Whitney class, since we work in cohomology $\mod 2$.
 
As in~\cite{kar2009mp}, we take large enough ball $B$ such that the section $f$ has no zeros on the set $SO(d)\times\partial B \subset E$. Now we can consider the obstruction as the relative Euler class of $f$ in the cohomology $H_{G'}^{2d-1}(SO(d)\times B, SO(d)\times \partial B)$.

Let us decompose $f$ into $s_U\oplus s_W$ by the corresponding $G'$-bundles. The section $s_U$ has no zeroes on $SO(d)\times \partial B$ and its Euler class resides in $e(s_U)\in H_{G'}^d(SO(d)\times B,SO(d)\times \partial B)$, the section $s_W$ has the Euler class $e(s_W)\in H_{G'}^{d-1} (SO(d)\times B) = H_{G'}^{d-1}(SO(d))$. The total Euler class is $e(f) = e(s_U)e(s_W)$, the product as defined as the natural product
$$
H_{G'}^d(SO(d)\times B,SO(d)\times \partial B)\times H_{G'}^{d-1} (SO(d)\times B)\to H_{G'}^{2d-1}(SO(d)\times B, SO(d)\times \partial B).
$$

Similar to~\cite{kar2009mp} (Lemma~1 in~\cite{kar2009mp} holds, its proof uses $G$-action without change), we obtain $e(s_U) = u\times 1\in H_{G'}^*(B\times SO(d), \partial B\times SO(d))$. Here we use the decomposition 
$$
H_{G'}^*(SO(d)\times B, SO(d)\times \partial B) = H^*(B,\partial B)\otimes H_{G'}^*(SO(d))
$$ 
and denote by $u$ the generator of $H^d(B, \partial B)$.

From the proof of Theorem~\ref{knaster-pow2} we know that $e(s_W)\not=0$, and by the K\"unneth formula and the multiplicative rule for the Euler class we have $e(f)=u\times e(s_W)\not=0$. The proof is complete.

\end{document}